\documentclass[a4paper,14pt]{article}

\usepackage{mathtools, nicematrix}
\usepackage[left=20mm, top=20mm, right=20mm, bottom=20mm, nohead]{geometry}
\usepackage[T1,T2A]{fontenc}
\usepackage[utf8]{inputenc}
\usepackage[english]{babel}
\usepackage{amsmath}
\usepackage{amsfonts,amssymb}
\usepackage{amsthm}
\usepackage[electronic]{ifsym}
\usepackage{authblk}
\usepackage{hyperref}
\usepackage[nameinlink, capitalize]{cleveref}
\usepackage[shortlabels]{enumitem}

\theoremstyle{theorem}
\newtheorem{theorem}{Theorem}[section]

\theoremstyle{definition}
\newtheorem{definition}[theorem]{Definition}
  
\theoremstyle{lemma}
\newtheorem{lemma}[theorem]{Lemma}

\theoremstyle{proposition}
\newtheorem{proposition}[theorem]{Proposition}

\theoremstyle{remark}
\newtheorem*{remark}{Remark}

\theoremstyle{corollary}
\newtheorem{corollary}[theorem]{Corollary}

\theoremstyle{conjecture}

\theoremstyle{definition}

\theoremstyle{question}
\newtheorem{question}[theorem]{Question}

\newcommand{\R}{\mathbb{R}}
\newcommand{\N}{\mathbb{N}}
\newcommand{\Q}{\mathbb{Q}}

\newcommand{\Z}{\mathbb{Z}}

\newcommand{\supp}{\mathrm{supp}}
\newcommand{\OR}{\mathrm{OR}}
\newcommand{\OSR}{\mathrm{OSR}}

\newcommand{\IP}{\mathrm{IP}}

\newcommand{\IOpen}{\mathsf{IOpen}}

\newcommand{\Exp}{\mathrm{Exp}}
\newcommand{\Log}{\mathrm{Log}}

\newcommand{\RCF}{\mathsf{RCF}}
\newcommand{\OF}{\mathsf{OF}}
\newcommand{\EF}{\mathsf{EF}}

\newcommand{\DOSR}{\mathsf{DOSR}}
\newcommand{\ZDOSR}{\mathbb{Z}\text{-}\mathsf{DOSR}}

\newcommand{\basis}{\mathrm{basis}}
\newcommand{\ad}{\mathrm{ad}}
\newcommand{\EndExt}{\subseteq^{\mathrm{end}}}
\newcommand{\teip}{\mathsf{TEIP}}

\begin{document}
\title{Some conservativity results for weak arithmetic theories in a language with exponentiation}
\author{Konstantin Kovalyov}
\affil{
Steklov Mathematical Institute of Russian Academy of Sciences, Moscow, Russia

\href{mailto:kovalyov.ka@mi-ras.ru}{kovalyov.ka@mi-ras.ru}
}

\maketitle

\begin{abstract}
Recall that $\IOpen$ denotes the theory of quantifier-free induction in the pure arithmetic language. In this paper, we consider an extension of $\IOpen$ in the arithmetical language expanded by the unary function $\exp: x \mapsto 2^x$, obtained by adding certain simple universal axioms for $\exp$, namely, 
$$\exp(1) = 2,  \quad \exp(x + y) = \exp(x) \exp(y), \quad \exp(x) > x.$$
We prove that this theory is $\Pi_1$-conservative over $\IOpen$ and $\forall$-conservative over the theory of discretely ordered $\Z$-semirings. This result extends a result of A. Wilkie \cite{wilkie1978} and provides a partial answer to a question of E. Jeřábek \cite{jerabek2024} concerning the strength of the theory of exponential integer parts. 
\end{abstract}

\section{Introduction}\label{section_introduction}

Recall that $\IOpen$ denotes the theory of discretely ordered semirings with quantifier-free induction in the pure arithmetic language. In this paper, we study an extension of $\IOpen$ in a language with the unary base-$2$ exponentiation $\exp$. Our main results are presented in \cref{theorem_on_extension} and \cref{corollary_on_extension}. We discuss the context and details of these results below.

The theory $\IOpen$ is well known and has been extensively studied; see, for instance, \cite{shepherdson:1964,wilkie1978, van_den_dries1980,van_den_dries1980_weak_arithm,otero1990, mourgues:1993, glivicka2017shepherdsons, kovalyov2024} and the references therein. In particular, $\IOpen$ is known for Shepherdson's Theorem \cite{shepherdson:1964}, which states that the models of $\IOpen$ are precisely the nonnegative integer parts of real closed fields (a discretely ordered subring $\mathcal M$ of an ordered ring $\mathcal R$ is called an \emph{integer part} if every element of $\mathcal R$ admits an integer part in $\mathcal M$). Moreover, M. H. Mourgues and J.-P. Ressayre \cite{mourgues:1993} proved that every real closed field admits an integer part. This yields an algebraic method for constructing models of $\IOpen$, which we will use later.

Several extensions of Shepherdson's Theorem have been studied; see, for instance, \cite{ITA_2008__42_1_105_0, jerabek2022, jerabek2024, kovalyov2026}. These papers employ the notion of an \emph{exponential integer part} of an \emph{exponential field} (see \cref{section_preliminaries} for the precise definitions). Note that every base-$2$ real closed exponential field (RCEF for short) satisfying the growth axiom admits an exponential integer part (see \cite{ressayre:1993}). In \cite{jerabek2024}, E. Jeřábek showed that the class of exponential integer parts of RCEFs satisfying the growth axiom is axiomatized by $\IOpen$ together with
\begin{align}
    \tag{$\Exp$}
    &\exp(1) = 2,  \quad \exp(x + y) = \exp(x) \exp(y), \quad \exp(x) > x, 
    \label{Exp_axioms}\\
    \tag{$\Log$}
    &(x > 0) \to \exists y \ \big(\exp(y) \leq x < \exp(y + 1) \big). \label{Log_axiom}
\end{align}
Following \cite{jerabek2024}, denote the theory $\IOpen + (\Exp) + (\Log)$ by $\teip_{2^x}^+$ .
This theory is not conservative over $\IOpen$: $\teip_{2^x}^+$ proves a purely arithmetical $\Pi_1$-sentence that is not provable in $\IOpen$. This naturally raises the question, posed in \cite{jerabek2024}, of how much stronger $\teip_{2^x}^+$ is than $\IOpen$. In particular, it is unknown whether $\teip_{2^x}^+$ is $\forall$-conservative over $\IOpen$.

In this paper, we consider the fragment $\IOpen + (\Exp)$ of $\teip_{2^x}^+$. Similar to \cite[Theorem 3.2]{jerabek2024}, it is easy to see that $\IOpen + (\Exp)$ is an axiomatization of the class of exponential integer parts of real closed \emph{weak exponential fields} satisfying the growth axiom (see \cref{subsection_exp_fields}). Hence, the theory $\IOpen + (\Exp)$ is quite natural.

We approach the above question from \cite{jerabek2024} by proving that the theory $\IOpen + (\Exp)$ is $\Pi_1$-conservative (and, hence, $\forall$-conservative) over $\IOpen$. The proof is model-theoretic: using Shepherdson's Theorem, we show that every model of $\IOpen$ admits an end-extension to a model of $\IOpen + (\Exp)$. Since A. Wilkie \cite{wilkie1978} proved that $\IOpen$ is $\forall$-conservative over the theory of discretely ordered $\Z$-semirings, the theory $\IOpen + (\Exp)$ is $\forall$-conservative over the theory of discretely ordered $\Z$-semirings. Hence, our result may also be viewed as an extension of Wilkie's theorem. Note that $\IOpen$ (and, hence, $\IOpen + (\Exp)$) is not $\Pi_1$-conservative over the theory of discretely ordered $\Z$-semirings, since, for instance, the sentence $$\forall x \ \exists y \leq x \ \big ( y^2 \leq x < (y + 1)^2 \big)$$ 
is not provable in the latter theory, but is provable in $\IOpen$.

\section{Preliminaries}\label{section_preliminaries}

\subsection{Conventions and notation}

By a ring, we mean an associative commutative ring with unity. Structures will usually be denoted by calligraphic letters (such as $\mathcal{M}, \mathcal{F}, \mathcal{R}, \dots$), and their domains by the corresponding letters $M, F, R, \dots$. Given a language $\mathcal L$, an $\mathcal L$-structure $\mathcal M$, and a function symbol $f \notin \mathcal L$, we denote by $\mathcal L(f)$ the expansion of $\mathcal L$ by $f$, and by $(\mathcal M, f_{\mathcal M})$ the expansion of $\mathcal M$ by the function $f_{\mathcal M}$. When no confusion can arise, we omit the subscript $\mathcal M$ from interpretations of symbols of $\mathcal L$ in the structure $\mathcal M$. For example, if
\[
\mathcal M = (M, +_\mathcal M, -_\mathcal M, \cdot_\mathcal M, 0_\mathcal M, 1_\mathcal M)
\]
is a ring, we may denote it by $(M, +, -, \cdot, 0, 1)$. 

We write $\mathcal M \subseteq \mathcal N$ to indicate that $\mathcal M$ is a substructure of $\mathcal N$. For a formula $\varphi$ and a tuple of variables $\vec x = (x_1, \dots, x_n)$, we write $\varphi(\vec x)$ to indicate that all free variables of $\varphi$ are among the variables in $\vec x$. An analogous convention applies to terms. In order to simplify notation, we will sometimes omit $\forall$-quantifiers at the beginning of formulas when writing axioms of a theory. For instance, we will write 
$x + y = y + x$ instead of $\forall x \ \forall y \ (x + y = y + x)$. We will write $T_1 + T_2$ for deductive closure of the union of theories $T_1$ and $T_2$. Given a theory $T$ in a language $\mathcal L$ and a theory $T' \supseteq T$ in a language $\mathcal L' \supseteq \mathcal L$, we say that $T'$ is \emph{$\forall$-conservative} over $T$ if for any universal $\mathcal L$-sentence $\varphi$, such that $T' \vdash \varphi$, we have $T \vdash \varphi$.

We will consider ordered semirings and ordered rings (the precise definitions will be given later in this section) in the languages
$$\mathcal L_{\OSR} := (+^{(2)}, \cdot^{(2)}, 0, 1, \leq^{(2)})$$
of ordered semirings and
$$\mathcal{L}_{\OR} := (+^{(2)}, -^{(1)}, \cdot^{(2)}, 0, 1, \leq^{(2)})$$
of ordered rings, respectively. We will denote the term
$$\underbrace{1 + (1 + (1 + \dots))}_{n \text{ times}}$$
simply by $n$. A formula $\varphi$ is said to be \emph{bounded}, or a \emph{$\Delta_0$-formula}, if every quantifier occurring in $\varphi$ is of the form $\exists x \leq t$ or $\forall x \leq t$, where $x$ does not occur in $t$. A formula $\varphi$ is said to be a \emph{$\Pi_1$-formula}, if it is of the form $\forall \vec x \ \theta$, where $\theta$ a $\Delta_0$-formula. 

If $\mathcal M = (M, \leq, \dots)$ is an ordered structure and $a \in M$, then we set
$$M^{>a} := \{m \in M \mid m > a\}.$$
For a subset $A \subseteq M$ and an element $a_0 \in M$, we write $A > a_0$ (respectively, $A \geq a_0$) if $a > a_0$ (respectively, $a \geq a_0$) for all $a \in A$. We write $A < B$ if $a < b$ for all $a \in A$ and $b \in B$.

If $\mathcal M = (M, +, -, \cdot, 0, 1, \leq)$ is an ordered ring, then we denote by $\mathcal M^+$ the semiring of nonnegative elements
$$(M^+, +, \cdot, 0, 1, \leq),$$
where
$$M^+ := \{a \in M \mid a \geq 0\}.$$

\subsection{Ordered rings and fields}

\begin{definition}
    $\mathsf{OR}$ is a theory in the language $\mathcal{L}_{\OR}$, consisting of the following axioms:
    \begin{enumerate}[start=0,label={(OR\arabic*)}]
        \item $x + (y + z) = (x + y) + z$;
        \item $x + y = y + x$;
        \item $x + 0 = x$;
        \item $x \cdot (y \cdot z) = (x \cdot y) \cdot z$;
        \item $x \cdot y = y \cdot x$;
        \item $x \cdot 1 = x$;
        \item $x \cdot (y + z) = x \cdot y + x \cdot z$;
        \item $x \leq x$;
        \item $(x \leq y \wedge y \leq x) \to x = y$;
        \item $(x \leq y \wedge y \leq z) \to x \leq z$;
        \item $x \leq y \vee y \leq x$;
        \item $x \leq y \to x + z \leq y + z$;
        \item $(x \leq y \wedge 0 \leq z) \to x \cdot z \leq y \cdot z$;
        \item $0 < 1$;
        \item $x + (-x) = 0$.
    \end{enumerate}
    Models of $\mathsf{OR}$ will be called \emph{ordered rings} (OR for short). 
\end{definition}

\begin{definition}
    $\mathsf{OF}$ is a theory in the language $\mathcal{L}_{\OR}$, consisting of the theory $\mathsf{OR}$ and the axiom
    $$(x \ne 0) \to \exists y \ (x \cdot y = 1).$$
    Models of $\mathsf{OF}$ will be called \emph{ordered fields} (OF for short). 
\end{definition}

\begin{remark}
    As usual, we will omit the symbol of multiplication $\cdot$ in the rest of the paper.
\end{remark}

\subsection{Discretely ordered rings and semirings}

\begin{definition}
    $\mathsf{DOR}$ is a theory in the language $\mathcal{L}_{\OR}$, consisting of the theory $\mathsf{OR}$ and the axiom
        $$x \leq 0 \vee 1 \leq x,$$
    which says that the order is discrete. Models of $\mathsf{DOR}$ will be called \emph{discretely ordered rings} (DOR for short). 
\end{definition}

\begin{definition}
    $\mathsf{DOSR}$ is a theory in the language $\mathcal{L}_{\OSR}$, consisting of the axioms (OR0)--(OR13) and the following axioms:
    $$x = 0 \vee 1 \leq x;$$
    $$x \leq y \leftrightarrow \exists z (x + z = y).$$
    Models of $\mathsf{DOSR}$ will be called \emph{discretely ordered semirings} (DOSR for short). 
\end{definition}

The theory $\mathsf{DOSR}$ will serve us as the basic arithmetic theory.

\begin{remark}
    One can freely move from models of $\mathsf{DOR}$ to $\mathsf{DOSR}$ and back. More precisely, for every DOR its nonnegative part is a DOSR and every DOSR is isomorphic to a nonnegative part of some DOR (just consider the ring of pairs $(m, n)$ modulo an equivalence relation $(m, n) \sim (m', n') :\iff m + n' = m' + n$). In the rest of the paper we are mostly dealing with DORs and their nonnegative parts without any mentioning of this equivalence.  For more details see, for instance, \cite{kaye}.
\end{remark} 

\begin{definition}
    $\ZDOSR$ is a theory in the language $\mathcal L_{\OSR}$, consisting of the theory $\DOSR$ and the axiom scheme
    $$\exists y \ \bigvee_{k = 0}^{n - 1} \big( x = n y + k\big), \quad n \in \N, n \geq 2.$$
    Models of $\ZDOSR$ will be called \emph{discretely ordered $\Z$-semirings} ($\Z$-DOSR for short). 
\end{definition}

\begin{definition}
    $\mathsf{IOpen}$ is a theory in the language $\mathcal{L}_{\OSR}$, consisting of the theory $\mathsf{DOSR}$ and the open induction schema
    \begin{align*}
        \Big(\varphi(0, \overline{y}) \wedge \forall x \big(\varphi(x, \overline{y}) \rightarrow \varphi(x + 1, \overline{y})\big)\Big) \rightarrow \forall x \ \varphi(x, \overline{y}),
    \end{align*}
    where $\varphi(x, \overline{y})$ is a quantifier-free formula in the language $\mathcal{L}_{\OSR}$.
\end{definition}

\begin{remark}
    It is easy to see that $\IOpen \vdash \ZDOSR$.
\end{remark}

\begin{theorem}[\cite{wilkie1978}]\label{wilkie_theorem}
    Every $\Z$-DOSR can be embedded into a model of $\IOpen$. Equivalently, $\IOpen$ is $\forall$-conservative over $\ZDOSR$.
\end{theorem}

\subsection{Real closed fields}

\begin{definition}
    $\mathsf{RCF}$ is a theory in the language $\mathcal{L}_{\OR}$, consisting of the theory $\mathsf{OF}$ and the scheme 
    $$(x < y \wedge t(x, \overline{w}) \leq 0 \wedge t(y, \overline{w}) \geq 0) \rightarrow \exists z ( (x \leq z \leq y) \wedge t(z, \overline{w}) = 0)$$
    for every $\mathcal L_{\OR}$-term $t(x, \overline{w})$. Models of $\mathsf{RCF}$ will be called \emph{real closed fields} (RCF for short).  
\end{definition}

For a detailed exposition of real closed fields, see, for instance, \cite{lang} and \cite{marker}. Here we state only the results that we will need.

\begin{definition}
    If $\mathcal F \subseteq \mathcal R$ are ordered fields, then $\mathcal R$ is called a \emph{real closure} of $\mathcal F$ if $\mathcal R$ is real closed and the extension $\mathcal F \subseteq \mathcal R$ is algebraic.
\end{definition}

\begin{theorem}[{\cite[Chapter XI, Theorem 2.9]{lang}}]\label{theorem_real_closure_exists_and_unique}
    For every ordered field $\mathcal F$ there exists its real closure $\mathcal R$. Moreover, this $\mathcal R$ is unique up to an isomorphism fixing $F$.
\end{theorem}

\begin{remark}
    Due to the theorem above one can speak about \emph{the} real closure of a given ordered field.
\end{remark}

\begin{lemma}[{\cite[Chapter XI, Lemma 3.4]{lang}}]\label[lemma]{theorem_relative_real_closure}
    If $\mathcal R$ is a real closed field and $\mathcal M \subseteq \mathcal R$ is an ordered subring, then the relative algebraic closure of $\mathcal M$ in $\mathcal R$ is the real closure of (the fraction field of) $\mathcal M$.
\end{lemma}

We now recall the following standard construction of formal series fields. Let $\mathcal F$ be an ordered field and $\mathcal P$ an ordered additive abelian group. Let $X$ be a formal variable. Denote by $X^\mathcal P$ the ordered multiplicative abelian group with domain $\{X^m \mid m \in P\}$ and operations defined in the obvious way. Thus, $X^\mathcal P$ is an ordered multiplicative copy of $\mathcal P$. For a function $f : P \to F$, denote by $\supp(f)$ the set
$$\{a \in P \mid f(a) \ne 0\},$$
and define
$$F((X^\mathcal P)) \ := \{f : P \to F \mid \text{the ordered set }\supp(f)\text{ is Noetherian}\}.$$
Elements of $F((X^\mathcal P))$ may be viewed as formal sums of the form
$$\sum_{m \in P} f(m) X^m.$$
One can define the field operations and the ordering on $F((X^\mathcal P))$ in a natural way as $f > 0$ iff $f\big(\max (\supp (f))\big) > 0$ (i.e. for every $m > 0$, $X^m$ is infinitely large); we denote the resulting ordered field by $\mathcal F((X^\mathcal P))$.

\begin{theorem}[{\cite[Theorem 4.3.7]{engler_prestel}}]\label{real_closed_fields_of_formal_sums}
    If an ordered field $\mathcal F$ is real closed and an ordered abelian additive group $\mathcal P$  is divisible, then $\mathcal F ((X^\mathcal P))$ is real closed. 
\end{theorem}

In \cref{section_conservativity}, this construction will be used to extend models of $\IOpen$.

We will also use the following construction to embed ordered abelian groups into ordered divisible abelian groups.

\begin{definition}
    Given an ordered additive abelian group $\mathcal P$, we denote by $\N^{-1}\mathcal P$ the \emph{ordered divisible hull} of $\mathcal P$, that is, the set $P \times \N^{>0}/\sim$, where $(m, k) \sim (m', k')$ if and only if $k'm = km'$, with the group operations and ordering defined in the natural way. The equivalence class of a pair $(m, k)$ is denoted by $\frac{m}{k}$.
\end{definition}

Clearly, $\N^{-1}\mathcal P$ is divisible.

\subsection{Integer parts and Shepherdson's Theorem}

\begin{definition}
    Let $\mathcal M$ and $\mathcal R$ be ordered rings, $\mathcal M \subseteq \mathcal R$ and $\mathcal M$ is DOR. Then $\mathcal M$ is called an \textit{integer part of} $\mathcal R$ if for all $r \in R$ there is $m \in M$ such that $m \leq r < m + 1$. Notation: $\mathcal M \subseteq^{\IP} \mathcal R$.
\end{definition}

\begin{remark}
    Since $\mathcal M$ is discretely ordered, for every $r \in R$ its integer part is uniquely defined.
\end{remark}

Let $\mathcal M$ be a discretely ordered ring. Denote by $\mathcal M^+$ the semiring of the nonnegative elements of $\mathcal M$, by $\mathcal{F}(\mathcal M)$ the quotient field of $\mathcal M$ and by $\mathcal{R}(\mathcal M)$ the real closure of $\mathcal{F}(\mathcal M)$ (which is unique up to an isomorphism).

\begin{theorem}[J. Shepherdson, \cite{shepherdson:1964}]\label{shepherdson_theorem}
    Let $\mathcal{M}$ be a discretely ordered ring. Then $\mathcal{M}^+ \vDash \mathsf{IOpen}$ iff $\mathcal{M} \subseteq^{\IP} \mathcal{R}(\mathcal{M})$.
\end{theorem}

According to \autoref{theorem_real_closure_exists_and_unique} and \cref{theorem_relative_real_closure}, we can reformulate \autoref{shepherdson_theorem} in the following form:

\begin{theorem}\label{shepherdson_theorem_reformulated}
    Let $\mathcal{M}$ be a discretely ordered ring. Then $\mathcal{M}^+ \vDash \mathsf{IOpen}$ iff there exists an ordered field $\mathcal{R}$ such that $\mathcal{M} \subseteq^{\IP} \mathcal{R}$ and $\mathcal{R} \vDash \RCF$.
\end{theorem}

\subsection{Extensions of $\mathsf{IOpen}$, exponential fields and exponential integer parts}\label{subsection_exp_fields}

\begin{definition}
    Denote by $(\Exp)$ the universal closure of the conjunction of the formulas
    $$\exp(1) = 2,  \quad \exp(x + y) = \exp(x) \exp(y), \quad \exp(x) > x,$$
    and by $(\Log)$ the universal closure of the formula
    $$(x > 0) \to \exists y \ \big(\exp(y) \leq x < \exp(y + 1) \big).$$
\end{definition}

\begin{remark}
    Here, $\exp$ stands for the unary base-2 exponentiation.
\end{remark}

We now explain the origin of these axioms.

\begin{definition}\label[definition]{def_exp_fields}
    $\EF^-$ is a theory in the language $\mathcal{L}_{\OR}(\exp)$, consisting of the theory $\mathsf{OF}$ and the following axioms:
    \begin{enumerate}[start=0,label={(EF\arabic*)}]
        \item $\exp(x) > 0$;
        \item $\exp(x + y) = \exp(x)\exp(y)$;
        \item $(x < y) \to (\exp(x) < \exp(y))$.
    \end{enumerate}
    $\EF$ is a theory in the language $\mathcal{L}_{\OR}(\exp)$, consisting of the theory $\mathsf{EF^-}$ and the axiom
    \begin{enumerate}[start=4,label={(EF\arabic*)}]
        \item $(x > 0) \to \exists y\ \big(x = \exp(y) \big)$.
    \end{enumerate}
    Models of $\EF^-$ will be called \emph{weak exponential fields} and models of $\EF$ will be called \emph{exponential fields}. 
\end{definition}

\begin{remark}
    It is not very hard to show that in fact (EF0) is provable in $\OF$ + (EF1)--(EF2). Some authors use different terminology: for instance, in \cite{vandendries_macintyre_marker:2001}, models of $\EF^-$ are called exponential fields.
\end{remark}

The standard examples of exponential fields are $\R_{\exp} := (\R, e^x)$ and $(\R, 2^x)$, where $\R$ denotes the ordered field of the real numbers.

\begin{definition}\label[definition]{def_exp_ip}
    Let $(\mathcal F, \exp)$ be a weak exponential field, $\mathcal M \subseteq \mathcal F$ be an ordered ring. Then $\mathcal M$ is called an \textit{exponential integer part} of $(\mathcal F, \exp)$ if $\mathcal M \subseteq^{\IP} \mathcal F$ and $M^+$ is closed under $\exp$ (i.e. for all $m \in M^+$ we have $\exp(m) \in M^+$). 
\end{definition}

\begin{theorem}[{\cite[Theorem 3.2]{jerabek2024}}]
    $\IOpen + (\Exp) + (\Log)$ is an axiomatization of the theory of exponential integer parts of models of $\RCF + \EF + \big(\exp(1) = 2\big) +  \forall x \ \big( \exp(x) > x\big)$. 
\end{theorem}

The following proposition is proved in the same way.

\begin{proposition}
    $\IOpen + (\Exp)$ is an axiomatization of the theory of exponential integer parts of models of $\RCF + \EF^- + \big(\exp(1) = 2\big) +  \forall x \ \big( \exp(x) > x\big)$. 
\end{proposition}

\begin{proof}
    The proof follows the same argument as that of \cite[Theorem 3.2]{jerabek2024}, based on Robinson's joint consistency theorem.
\end{proof}

\subsection{End-extensions of discretely ordered semirings}

\begin{definition}
    Let $\mathcal M$ and $\mathcal N$ be discretely ordered semirings. Then $\mathcal N$ is called an \emph{end-extension} of $\mathcal M$ if $\mathcal M \subseteq \mathcal N$ and for every $m \in M$ and $n \in N$, $n \leq m$ implies $n \in M$. Notation: $\mathcal M \EndExt \mathcal N$.

    For discretely ordered rings $\mathcal M$ and $\mathcal N$, we say that $\mathcal N$ is an \emph{end-extension} of $\mathcal M$, and write $\mathcal M \EndExt \mathcal N$, if $\mathcal M^+ \EndExt \mathcal N^+$.
\end{definition}

The following facts are standard.

\begin{proposition}[{\cite[Theorem 2.7]{kaye}}]
    Let $\mathcal M$ and $\mathcal N$ be discretely ordered semirings with $\mathcal M \EndExt \mathcal N$. Then for every $\Delta_0$-formula $\theta(\vec x)$ and every $\vec a \in M$,
    $$\mathcal M \vDash \theta(\vec a) \ \iff \ \mathcal N \vDash \theta(\vec a).$$
\end{proposition}

\begin{corollary}
    Let $T$ be a theory in the language $\mathcal L_{\OSR}$ and let $U$ be a theory in an expansion of the language $\mathcal L_{\OSR}$ such that $\DOSR \subseteq T, U$. Suppose that every model of $T$ admits an end-extension to a model of $U$. Then $U$ is $\Pi_1$-conservative over $T$, that is, every $\Pi_1$-sentence provable in $U$ is also provable in $T$.    
\end{corollary}

\section{Conservativity results for $\IOpen + (\Exp)$}\label{section_conservativity}

The main result of this section is as follows. 

\begin{theorem}\label{theorem_on_extension}
    Every model of $\IOpen$ admits an end-extension to a model of $\IOpen + (\Exp)$. Consequently, the theory $\IOpen + (\Exp)$ is $\Pi_1$-conservative over $\IOpen$. 
\end{theorem}

Before the proof, we observe the following known construction, which allows us to extend models of $\IOpen$.  

\begin{proposition}[see \cite{shepherdson:1964, mourgues:1993}]\label[proposition]{prop_on_ext_of_models_iopen}
    Let $\mathcal M$ be an ordered ring such that $\mathcal M^+ \vDash \IOpen$, and let $\mathcal P$ be an ordered additive abelian divisible group. Denote by $\mathcal R$ the real closure of the fraction field of $\mathcal M$, and by $\mathcal N$ the ordered ring with domain
    $$N \ := \ \{f \in \mathcal R((X^\mathcal P)) \mid \supp(f) \geq 0, \ f(0) \in M\}.$$
    Then $\mathcal M \EndExt \mathcal N$ and $\mathcal N^+ \vDash \IOpen$.
\end{proposition}

\begin{proof}
    By \cref{real_closed_fields_of_formal_sums}, the ordered field $\mathcal R((X^\mathcal P))$ is real closed since $\mathcal P$ is divisible. Clearly, $\mathcal N$ is discretely ordered and $\mathcal M \EndExt \mathcal N$. Since $\mathcal N$ is an integer part of the real closed field $\mathcal R((X^\mathcal P))$, by \cref{shepherdson_theorem_reformulated} we obtain $\mathcal N^+ \vDash \IOpen$.
\end{proof}

\begin{proof}[Proof of \cref{theorem_on_extension}]
    Let $\mathcal M$ be a discretely ordered ring such that $\mathcal M^+ \vDash \IOpen$. We construct recursively a sequence of models 
    $$\mathcal N_0 \subseteq \mathcal N_1 \subseteq \mathcal N_2 \subseteq \dots,$$
    sets $P_{i + 1} \subseteq \mathcal N_{i + 1}$ and maps $\exp_i : N^+_i \to N^+_{i + 1}$ for all $i \in \N$ such that for all $i \in \N$ the following conditions hold:
    \begin{enumerate}[start=0,label={\rm (\arabic*)}]
        \item $\mathcal N_0 = \mathcal M$;
        \item $\mathcal N_{i + 1}^+ \vDash \IOpen$ and $\mathcal N_i \EndExt \mathcal N_{i + 1}$;
        \item $(P_{i + 1}, +, 0, \leq)$ is a subgroup of $(N_{i + 1}, +, 0, \leq)$;
        \item $\mathcal N_{i + 1} = \mathcal N_i \oplus P_{i + 1}$ as abelian groups; moreover, for all $m_0 \in N_i$ and $m_1 \in P_{i + 1}^{>0}$, we have $m_0 < m_1$;
        \item $\exp_i$ satisfies the axiom $(\Exp)$ on $N_i^+$;
        \item $\exp_i|_{N_{i - 1}^+} = \exp_{i - 1}$ for $i > 0$.
    \end{enumerate}

    First, let us describe the idea informally. On every step, the model $\mathcal N_i$ is decomposed into a ``small'' part ($\mathcal N_{i - 1}$, if $i > 0$, and $\Z$, if $i = 0$) and a ``big'' part ($P_i$, if $i > 0$ and $ N^{\basis}_0$, if $i = 0$). We already have the exponential function constructed on the nonnegative elements of the ``small'' part of $N_i$. Then we construct a model $\mathcal N_{i + 1}$ from $\mathcal N_i$ as in \cref{prop_on_ext_of_models_iopen} and decompose it into its ``small'' part $\mathcal N_i$ and its ``big'' part $P_{i + 1}$. After that, we extend the exponential function on $\mathcal N_i^+$ by sending an element $m$ of the ``big'' part of $N_i$ to the element $X^m \in P_{i + 1}$ of the ``big'' part of $\mathcal N_{i + 1}$.

    Now let us proceed with the construction.  $\mathcal N_0$ is already constructed. On the $i$-th step we construct $\mathcal N_{i + 1}$, $P_{i + 1}$ and $\exp_i : N_i^+ \to N_{i + 1}^+$.

    We begin with the base case $i = 0$. Denote by $\mathcal R_0$ the real closure of the fraction field of $\mathcal N_0$, and denote by $\mathcal N_0^{\ad}$ the additive ordered group $(N_0, +, 0, \leq)$. By \cref{shepherdson_theorem}, $\mathcal N_0 \subseteq^{\IP} \mathcal R_0$. Let $N_0^{\basis}$ be a maximal subset of $N_0^+ \setminus \{1\}$ such that $\{1\} \cup N_0^{\basis}$ is linearly independent over $\Q$ in $\N^{-1} \mathcal N_0^{\ad}$; such exists by Zorn's Lemma. Clearly, every element of $N_0$ has a unique representation as a linear combination of elements of $\{1\} \cup N_0^{\basis}$ with rational coefficients (though not every such rational linear combination need lie in $N_0$). Now define a discretely ordered ring $\mathcal N_1$ with domain
    $$N_1\ :=\ \{f \in R_0((X_0^{\N^{-1} \mathcal N_0^{\ad}})) \mid \supp (f) \geq 0 \text{ and } f(0) \in N_0\},$$
    and define
    $$P_1 := \{f \in R_0((X_0^{\N^{-1} \mathcal N_0^{\ad}})) \mid \supp (f) > 0\}.$$
    Clearly, the conditions (2) and (3) hold. By \cref{prop_on_ext_of_models_iopen}, $\mathcal N_1^+ \vDash \IOpen$ and $\mathcal N_0 \EndExt \mathcal N_1$; hence, condition (1) also holds.  It remains to define $\exp_0 : N_0^+ \to N_1^+$. Let $m \in N_0^+$. Then $m$ is uniquely represented as $m = q_0 + \sum\limits_{j = 1}^n q_j m_j$ with $q_0, \dots, q_n \in \Q$, $m_1, \dots, m_n \in N_0^{\basis}$, and $n \in \N$. Note that $\sum\limits_{j = 1}^n q_j m_j$ is either zero or nonstandard (i.e., larger than any standard integer) and hence nonnegative. Now define 
    $$\exp_0(m) := 2^{q_0} X_0^{\sum\limits_{j = 1}^n q_j m_j}.$$
    Since $2^{q_0}$ is real algebraic, it belongs to $R_0$, and so $\exp_0(m) \in N_1^+$. One can readily verify that $\exp_0$ satisfies $(\Exp)$ on its domain, i.e., condition (4) holds.

    Now consider the inductive step. The construction is analogous to the base case. Suppose, as the induction hypothesis, that $\mathcal N_{i + 1}$, $P_{i + 1}$, and $\exp_i : N_i^+ \to N_{i + 1}^+$ have been constructed satisfying conditions (1)–(5). Denote by $\mathcal R_{i + 1}$ the real closure of the fraction field of $\mathcal N_{i + 1}$. Now define a discretely ordered ring $\mathcal N_{i + 2}$ with domain
    $$N_{i + 2}\ :=\ \{f \in R_{i + 1}((X_{i + 1}^{\N^{-1} \mathcal P_{i + 1}})) \mid \supp f \geq 0 \text{ and } f(0) \in N_{i + 1}\},$$
    where $\mathcal P_{i + 1} = (P_{i + 1}, +, 0, \leq)$, and define
    $$P_{i + 2} := \{f \in R_{i + 1}((X_{i + 1}^{\N^{-1} \mathcal P_{i + 1}})) \mid \supp f > 0\}.$$
    Clearly, the conditions (2) and (3) hold. By \cref{prop_on_ext_of_models_iopen}, $\mathcal N_{i + 2}^+ \vDash \IOpen$ and $\mathcal N_{i + 1} \EndExt \mathcal N_{i + 2}$, and so condition (1) also holds. It remains to define $\exp_{i + 1} : N_{i + 1}^+ \to N_{i + 2}^+$. Extend $\exp_i$ to $N_i$ by $\exp_i(-m) := \frac{1}{\exp_i(m)} \in R_{i + 1}$ for $m \in N_i$. Then for $m \in N_{i + 1}^+$ define $\exp_{i + 1}(m)$ as $\exp_i(m_0) X_{i + 1}^{m_1}$, where $m_0 \in N_i$ and $m_1 \in P_{i + 1}$  are the unique elements such that $m = m_0 + m_1$ (note that $m_1 \geq 0$ by condition (3)). Clearly, the conditions (4) and (5) are also satisfied. 
    
    Finally, set $\mathcal N := \bigcup\limits_{i \in \N} \mathcal N_i$ and $\exp := \bigcup\limits_{i \in \N} \exp_i : N^+ \to N^+$. Then $\mathcal N^+ \vDash \IOpen$ since $\IOpen$ is $\forall\exists$-axiomatizable and unions of chains preserve $\forall\exists$-sentences (see, e.g., \cite[Lemma 3.1.6]{tent2012course}), and $(\mathcal N^+, \exp) \vDash (\Exp)$ since $(\Exp)$ is purely universal. Clearly, $\mathcal M \EndExt \mathcal N$.
\end{proof}

\begin{corollary}\label[corollary]{corollary_on_extension}
    The theory $\IOpen + (\Exp)$ is $\forall$-conservative over the theory $\ZDOSR$. 
\end{corollary}

\begin{proof}
    Apply \cref{theorem_on_extension} and \cref{wilkie_theorem}.
\end{proof}

\begin{corollary}
    The formula $(x > 0) \to \exists y \ \big(\exp(y) \leq x < \exp(y + 1) \big)$ is not provable in $\IOpen + (\Exp)$. 
\end{corollary}

\begin{proof}
    Suppose that $\IOpen + (\Exp) \vdash (x > 0) \to \exists y \ \big(\exp(y) \leq x < \exp(y + 1) \big)$. This implies that $\IOpen + (\Exp)$ and $\teip_{2^x}^+$ are equivalent. But by \cref{theorem_on_extension}, $\IOpen + (\Exp)$ is $\Pi_1$-conservative over $\IOpen$ while $\teip_{2^x}^+$ is not (see \cite[Example 4.9]{jerabek2024}).
\end{proof}

However, the question of the $\forall$-conservativity of $\teip_{2^x}^+$ over $\IOpen$, posed in \cite{jerabek2024}, remains unanswered. The construction used in the proof of \cref{theorem_on_extension} does not allow us to add integer parts of logarithms to models of $\IOpen$. We leave this question and the following related questions open.

\begin{question}
    \begin{enumerate}[(1)]
        \item Is  $\IOpen + (\Exp)$ conservative over $\IOpen$?
        \item More generally, what is the $\mathcal L_{OSR}$-fragment of $\IOpen + (\Exp)$?
        \item  Is $\teip_{2^x}^+$ $\forall$-conservative over $\IOpen + (\Exp)$ in the language $\mathcal L_{\OSR}(\exp)$?
    \end{enumerate}
\end{question}

\section*{Funding}
This work was supported by the Russian Science Foundation under grant no.~25-21-00632; see \href{https://rscf.ru/project/25-21-00632/}{https://rscf.ru/project/25-21-00632/}.

\bibliographystyle{unsrt}
\bibliography{references}

\end{document}